\documentclass[11pt]{article}

\usepackage{epic,latexsym,amssymb}
\usepackage{color}
\usepackage{tikz}
\usepackage{amsfonts,epsf,amsmath}
\usepackage{epic,latexsym,amssymb,xcolor}
\usepackage{xcolor}
\usepackage{tikz}
\usepackage{amsfonts,epsf,amsmath,leftidx}
\usepackage{float}
\usepackage{mathrsfs}
\usetikzlibrary{arrows}

\usepackage[bottom=1.5in, top=1.2in, left=1.5in, right=1.5in]{geometry}

\newtheorem{them}{Theorem}
\newtheorem{lem}{Lemma}
\newtheorem{quest}{Question}

\newtheorem{defn}{Definition}
\newtheorem{cl}{Claim}
\newtheorem{cor}{Corollary}

\newtheorem{prop}{Proposition}
\newtheorem{conj}{Conjecture}

\newtheorem{problem}{Problem}

\newcommand{\qed}{$\Box$}

\newcommand{\smallqed}{{\tiny ($\Box$)}}

\newcommand{\proof}{\noindent\textbf{Proof. }}

\let\oldenumerate\enumerate
\renewcommand{\enumerate}{
  \oldenumerate
  \setlength{\itemsep}{0pt}
  \setlength{\parskip}{0pt}
  \setlength{\parsep}{0pt}
}

\begin{document}

\title{On graphs whose domination number is equal to chromatic and dominator chromatic numbers}

\author{$^{\dagger}$David A. Kalarkop, $^{\ddagger, \S}$Pawaton Kaemawichanurat  and $^{\dagger}$Raghavachar Rangarajan
\\ \\
$^{\dagger}$Department of Studies in Mathematics,\\ University of Mysore,\\ Manasagangothri, Mysuru-570006, India.\\
$^{\ddagger}$Department of Mathematics, Faculty of Science\\
King Mongkut's University of Technology Thonburi,\\
Thailand\\
$^{\S}$Mathematics and Statistics with Application (MaSA)
\\
\small \tt Email: david.ak123@gmail.com, pawaton.kae@kmutt.ac.th, rajra63@gmail.com }

\date{}
\maketitle


%

\begin{abstract}
For a graph $G = (V(G), E(G))$, a dominating set $D$ is a vertex subset of $V(G)$ in which every vertex of $V(G) \setminus D$ is adjacent to a vertex in $D$. The domination number of $G$ is the minimum cardinality of a dominating set of $G$ and is denoted by $\gamma(G)$. A coloring of $G$ is a partition $C = (V_{1}, ... ,V_{k})$ such that each of $V_{i}$ in an independent set. The chromatic number is the smallest $k$ among all colorings $C = (V_{1}, ... ,V_{k})$ of $G$ and is denoted by $\chi(G)$. A coloring $C = (V_{1}, ... ,V_{k})$ is said to be dominator if, for all $V_{i}$, every vertex $v \in V_{i}$ is singleton in $V_{i}$ or is adjacent to every vertex of $V_{j}$. The dominator chromatic number of $G$ is the minimum $k$ of all dominator colorings of $G$ and is denoted by $\chi_{d}(G)$. Further, a graph $G$ is $D(k)$ if $\gamma(G) = \chi(G) = \chi_{d}(G) = k$. In this paper, for $n \geq 4k - 3$, we prove that there always exists a $D(k)$ graph of order $n$. We further prove that there is no planar $D(k)$ graph when $k \in \{3, 4\}$. Namely, we prove that, for a non-trivial planar graph $G$, the graph $G$ is $D(k)$ if and only if $G$ is $K_{2, q}$ where $q \geq 2$.
\end{abstract}

{\small \textbf{Keywords:} domination number, chromatic number, dominator chromatic number, planar graphs.} \\
\indent {\small \textbf{AMS subject classification: 05C15, 05C69, 05C10} }

\section{Introduction}
By a graph $G$, we mean a simple connected graph where $V(G)$ is the vertex set and $E(G)$ is the edge set of $G$. If two vertices $v_1$ and $v_2$ are adjacent in $G$, then we say $e=v_1 v_2$ is an edge and we write $v_1 v_2 \in E(G)$. If $S \subset V(G) $, then the graph induced on $S$, denoted by $G[S]$ is the graph with vertex set $S$ and edge set $E=\{e=u v \in E(G): u,v \in S \}$. A graph is planar if it can be embedded in the plane. For more information on graph theoretic definitions and terminologies, we refer to \cite{CD}.

Coloring of a graph and domination in graphs are two important and fascinating areas in graph theory due to its wide applications in other areas such as computer science, communication networks etc. A coloring of a graph is the assignment of colors to the vertices of $G$ such that no two adjacent vertices receive the same color. The minimum number of colors used among all the colorings of $G$ is said to be \textit{chromatic number of $G$}, denoted by $\chi(G)$. A subset $D$ of $V(G)$ is said to be a \textit{dominating set} of $G$ if every vertex of $G$ is in $D$ or has a neighbor in $D$. The minimum cardinality of a dominating set of $G$ is said to be \textit{domination number} of $G$, denoted by $\gamma(G)$. A dominating set $D$ of $G$ such that the subgraph induced by $D$ has no isolate vertices is said to be a \textit{total dominating set} of $G$. The minimum cardinality of a total dominating set of $G$ is said to be \textit{total domination number} of $G$, denoted by $\gamma_{t}(G)$. Let $D$ be the dominating set of $G$ and let $u \in D$, a vertex $v$ is said to be an \emph{external private neighbor} of $u$ with respect to $D$ if $v$ is adjacent to only $u$ in $D$. For more details on domination and coloring of graphs, refer \cite{AB,IJ,KL}.

A graph is said to be $k$\emph{-colorable} if it can be colored using $k$-colors. A $k$\emph{-coloring} $C=(V_1, V_2, \cdots, V_k)$ of a $k$-colorable graph partitions $V(G)$ into independent sets $V_i$, for $1 \leq i \leq k$ and each $V_i$ is said to be a color class with respect to coloring $C$. A vertex $u \in V(G)$ is said to \textit{dominate the color class $V_i$} (or $u$ is said to be a \textit{dominating vertex of the color class $V_i$}) if $u$ is adjacent to all the vertices of $V_i$. There are several coloring problems that are defined involving the concept of domination, two such coloring problems are \textit{dominator coloring} and \textit{dominated coloring} of a graph.  Dominator coloring of a graph was studied for the first time by Gera et al.\cite{GH}. Dominator coloring of $G$ is a coloring $C$ of $G$ such that every vertex dominates at least one color class in $C$. The minimum number of colors used among all the dominator colorings of $G$ is said to be \textit{dominator chromatic number} of $G$, denoted by $\chi_{d}(G)$.
Dominated coloring of a graph was introduced by Merouane et al.\cite{MN}. Dominated coloring of $G$ is a coloring in which every color class has a dominating vertex. The minimum number of colors used among all dominated colorings of $G$ is said to be \textit{dominated chromatic number} of $G$, denoted by $\chi_{dom}(G)$.
\vskip 5 pt

\indent We further define the following graph,
\begin{defn}
	A graph $G$ such that $\chi(G)=\gamma(G)=\chi_{d}(G)=k$ is said to be $D(k)$ graph.
\end{defn}
\vskip 5 pt

\indent For the study on $D(k)$ graphs, when $k = 1$, the only one $D(1)$ graph is the singleton vertex. For the case when $k=2$, Gera \cite{EF} has proved the following result,
\begin{prop}\label{d2}\cite{EF}
	Let $G$ be a connected graph of order $n$. Then
	$\chi_{d}(G) = 2$ if and only if $G = K_{a,b}$ for $a, b \in \mathbb{N}$.
\end{prop}

In the same paper, Gera posed the following question,
\begin{quest}
	For what graphs does satisfy $\chi(G)=\gamma(G)=\chi_{d}(G)=k$ ?
\end{quest}
\vskip 5 pt

\noindent For related results in the study of domination and chromatic numbers of a graph, Arumugam et al. \cite{12}, Chellali et al. \cite{34} have given the characterization of some classes of graphs with specific domination number and chromatic number. Recently, Goddard et al.\cite{56} have proved some results on dominator colorings in planar graphs with small diameter.  In this paper, we initiate the study on graphs whose domination number is equal to chromatic and dominator chromatic numbers.

\section{Main Results}
In this section, we present the realizability results on $D(k)$ graphs of order $n \geq 4k - 3$. We establish the characterization of $D(3)$ graphs and that there are no planar graphs which are $D(k)$, for $k \geq 3$. All the main results in this section will be proved in Section \ref{proofs}. Our first main theorem is useful throughout all the proofs in this paper.
\vskip 5 pt

\begin{them}\label{thm 1}
If $G$ is a $D(k)$ graph with any $\chi_{d}$-coloring $C$, then every color class in $C$ is dominated by a vertex and every vertex dominates exactly one color class in $C$.
\end{them}

\indent By Theorem \ref{thm 1}, if $G$ is a $D(k)$ graph, then any $\chi_{d}$-coloring of $G$ is a $\chi_{dom}$-coloring of $G$. Therefore $\chi_{dom}(G) \leq \chi_{d}(G)$. Also $\chi(G) \leq \chi_{dom}(G) \leq \chi_{d}(G)=\chi(G)$ implying that $\chi_{d}(G)=\chi_{dom}(G)$. Thus, we have Corollary \ref{cor 1} as follows.
\vskip 5 pt

\begin{cor}\label{cor 1}
If $G$ is a $D(k)$ graph. Then  $\gamma(G)=\chi(G)=\chi_{d}(G)=\chi_{dom}(G)=k$.	
\end{cor}
\vskip 5 pt

\indent Since a total dominating set is a dominating set, $\gamma(G) \leq \gamma_{t}(G)$. Thus, if $G$ has a total dominating set of size $\gamma(G)$, then $\gamma_{t}(G) = \gamma(G)$. Hence, from the proof of Corollary \ref{cor 1}, we obtain the following corollary.
\vskip 5 pt

\begin{cor}
If $G$ is a $D(k)$ graph with $\chi_{d}$-coloring $C=(V_1, V_2, \cdots, V_k)$, then there exists a set $\{x_i: x_i \in V_i\}$ which is a total dominating set of $G$. In particular, if $G$ is a $D(k)$ graph, then
$$\gamma(G)=\gamma_{t}(G)=\chi(G)=\chi_{d}(G)=\chi_{dom}(G)=k.$$
\end{cor}

\indent Our next main result is to establish the existence of $D(k)$ graphs of order $n \geq 4k - 3$.

\noindent Let
\vskip 5 pt

\indent $\mathcal{D}(k, n):$ the class of $D(k)$ graphs of order $n$.
\vskip 5 pt

\begin{them}\label{real}
For any integer $k \geq 2$ and natural number $n \geq 4k - 3$, we have that $\mathcal{D}(k, n) \neq \emptyset$.
\end{them}

\indent Next, we will establish the construction of $D(3)$ graphs.
\vskip 5 pt

\noindent \textbf{The class} $\mathcal{D}(3)$
\vskip 5 pt

\indent Let $V_1, V_2, V_3$ be three independent sets with $|V_1|\geq 3, |V_2|\geq 3$ and $|V_3|=1$. For $1 \leq i \leq 3$, let $x_i \in V_i $ and let vertices $y_1,y_2,y_3$ be such that $x_i$ is adjacent to $y_j$ if and only if $i = j$. Further, we let $y_1, y_3 \in V_2$, $y_2\in V_1$. A graph in the class $\mathcal{D}(3)$ is obtained from $V_1, V_2, V_3$ by adding edges as follows.
\begin{enumerate}
  \item Join $y_1$ to every vertex in $V_1$ and join $y_2$ to every vertex in $V_2$.
  \item Join $x_1$ to $x_3$ and join $x_{1}$ to some (or all) the vertices of $V_{2}-\{y_3\}$.
  \item Join $y_3$ to some (or all) the vertices of $V_{1}-\{x_1\}$.
  \item Join each vertex in $V_{1}-\{x_1,y_2\}$ to every vertex in $V_2$ or $V_{3}$ (but not both) and join each vertex in $V_{2}-\{y_3, y_1\}$ to every vertex in $V_2$ or $V_{3}$ (but not both) in such a way that $x_3$ has at least two non-neighbors in both $V_1$ and $V_2$.
  \item For all $x \in V_1$, $y \in V_2$ such that $xy \notin E(G)$, either $x$ has at least one non-neighbor other than $y$ in $V_2$ or $y$ has at least one non-neighbor other than $x$ in $V_1$.
\end{enumerate}

\noindent We prove that:

\begin{them}\label{thm d3}
Let $G$ be a graph. Then, $G$ is a $D(3)$ graph if and only if $G \in \mathcal{D}(3)$. In particular, there is no planar $D(3)$ graph.
\end{them}

\noindent Further, we prove that there is no $D(4)$ graph $G$ when the graph $G$ is planar.

\begin{them}\label{thm 2}
There is no planar $D(4)$ graph.
\end{them}

\indent For a planar graph $G$, we have that $\chi(G)\leq 4$. When $k=1$, $G$ is the trivial graph. When $k=2$, Proposition \ref{d2} implies that $G$ is $K_{p,q}$, where $p,q \in \mathbb{Z^+}$. If $p=1$, then $\gamma(G)=1$, a contradiction. If $p,q\geq 3$, then $K_{p,q}$ has a minor of $K_{3,3}$ as a subgraph, contradicting $G$ is planar. Therefore, for $p=2$ and $q \geq 2$. By Theorems \ref{thm d3} and \ref{thm 2}, we obtain the following theorem.
\vskip 5 pt

\begin{them}\label{thm 5}
Let $G$ be a non-trivial planar graph. Then, $G$ is $D(k)$ if and only if $k = 2$ and $G\cong K_{2,q}$ when $q \geq 2$.
\end{them}

\section{Proofs}\label{proofs}

\subsection{Proof of Theorem \ref{thm 1}} Let $G$ be a $D(k)$ graph and $C=(V_1, V_2, \cdots, V_k)$ be a $\chi_{d}$-coloring of $G$. Since $\chi(G)=k$, there is a vertex $x_1$ in $V_1$ that has at least one neighbor, say $x_j$ in each $V_j$, for all $2 \leq j \leq k$. Now for $1 \leq i \leq k$, the set $D=\{x_i: x_{i}\in V_i\}$ is a minimum total dominating set of $G$ and so every vertex of $D$ has an external private neighbor (say $y_i$ is the external private neighbor of $x_i$). Then $y_i$ dominates $V_i$. Therefore, every color class is dominated. Now suppose there exists a vertex $z$ which dominates at least two color classes in $C$. Let $z \in V_1$ such that $z$ dominates $V_2$ and $V_3$. Then $D^{\prime}=\{z, y_i: 1\leq i \leq k$ and $i \neq 2,3\}$ is a dominating set of $G$ with cardinality $k-1$, a contradiction. Therefore every vertex dominates exactly one color class in $C$. This proves Theorem \ref{thm 1}.

\subsection{Proof of Theorem \ref{real}}
\indent We prove by constructions. For $k=1$, $D(1)$ is the trivial graph. For $k=2$, a complete bipartite graph $K_{p, q}$ where $p, q \geq 2$ is the $D(2)$ graph. Then, we may suppose that $k \geq 3$. In the following, we establish realizability of the existence of $D(k)$ graphs for all $k \geq 3$ with given order. We may establish constructions of $D(k)$ graphs depending on parity of $k$.
\vskip 5 pt

\noindent \textbf{The graph} $D^{odd}_{k, n}$
\vskip 5 pt

\noindent For an odd integer $k \geq 3$, we let
\vskip 5 pt

\indent $X=\{x_1, x_2, \cdots, x_{k-1} \},$
\vskip 5 pt

\indent $Y=\{y_1, y_2, \cdots, y_{k-1} \},$
\vskip 5 pt

\indent $Z=\{z_1, z_2, \cdots, z_{k-1} \},$
\vskip 5 pt

\indent $W=\{w_1, w_2, \cdots, w_{k-1} \}$
\vskip 5 pt

\noindent and, for integer $t \geq 1$, we let
\vskip 5 pt

\begin{equation*}
U = \left\{
\begin{array}{rl}
	\emptyset & \text{if } n = 4k - 3,\\\\
	\{u_{1}, ..., u_{t}\} & \text{if } n = 4k - 3 + t.
\end{array} \right.
\end{equation*}

\noindent Then, for $2 \leq i \leq k-1$, we let $P_i=\{x_i, y_i, z_i,w_i\}$ and $P_{1} = \{x_1, y_1, z_1,w_1\} \cup U$ be independent sets and we let $P_k=\{x_k\}$. The graph $\mathcal{D}^{odd}_{k, n}$ is obtained from $P_{1}, ..., P_{k}$ by adding edges as follows:
\vskip 5 pt
\begin{enumerate}
  \item For all $1 \leq i \leq \frac{k-1}{2}$, join every vertex in $\{x_{2i},y_{2i},z_{2i}\}$ to every vertex in $\{x_{2i-1},y_{2i-1},z_{2i-1}\}$. Further, join $w_{2i-1}$ to $y_{2i},z_{2i}$ and $w_{2i}$ to $y_{2i-1},z_{2i-1}$.
  \item Add edges so that $X\cup \{x_{k}\}$ is a $k$-clique.
  \item Join $x_k$ to all the vertices of $W \cup U$.
  \item If $U \neq \emptyset$, then join every vertex in $U$ to $y_{2}$ and $z_{2}$.
\end{enumerate}

\noindent Observe that, by (a), the subgraph induced by $P_{2i-1} \cup P_{2i}$ is $K_{3, 3}$ for all $1 \leq i \leq \frac{k-1}{2}$. An example of the graph $D^{odd}_{3, 9}$. is illustrated by Figure \ref{g39}.

\begin{figure}[H]
		\begin{center}
		\unitlength 1mm 
		\linethickness{0.4pt}
		\ifx\plotpoint\undefined\newsavebox{\plotpoint}\fi 
		\begin{picture}(92.375,101.75)(0,0)
			\put(17.5,75.5){\circle*{2.693}}
			\put(17.75,60.75){\circle*{2.693}}
			\put(17.5,46.75){\circle*{2.693}}
			\put(34.75,47){\circle*{2.693}}
			\put(34.25,61.75){\circle*{2.693}}
			\put(33.75,76.25){\circle*{2.693}}
			\put(53.75,76.25){\circle*{2.693}}
			\put(36.25,29.5){\circle*{2.693}}
			\put(16.75,29.5){\circle*{2.693}}
			\qbezier(16.75,75.5)(33.25,101.75)(53.75,77)
			\put(53.75,77){\line(-1,0){20.25}}
			\put(33.5,77){\line(1,0){.5}}
			\multiput(53.75,76.25)(-.0336914063,-.0913085938){512}{\line(0,-1){.0913085938}}
			\multiput(17.75,76.25)(.0391566265,-.0337349398){415}{\line(1,0){.0391566265}}
			\multiput(17.25,75.75)(.0337301587,-.0560515873){504}{\line(0,-1){.0560515873}}
			\multiput(18.5,60.75)(.0336757991,.0371004566){438}{\line(0,1){.0371004566}}
			\multiput(17.25,60.75)(2.09375,.03125){8}{\line(1,0){2.09375}}
			\multiput(17,60.75)(.0430174564,-.0336658354){401}{\line(1,0){.0430174564}}
			\put(17.5,48){\line(1,0){16.25}}
			\multiput(15.5,46.75)(.0409292035,.0337389381){452}{\line(1,0){.0409292035}}
			\multiput(34.25,76.75)(-.0337186898,-.0568400771){519}{\line(0,-1){.0568400771}}
			\put(14.25,77){\makebox(0,0)[cc]{$x_1$}}
			\put(31.5,79){\makebox(0,0)[cc]{$x_2$}}
			\put(57.25,80){\makebox(0,0)[cc]{$x_3$}}
			\put(14.5,44.25){\makebox(0,0)[cc]{$z_1$}}
			\put(35.5,44.25){\makebox(0,0)[cc]{$z_2$}}
			\put(36,30){\framebox(0,.25)[]{}}
			\multiput(35.75,29.75)(-.0361271676,.0337186898){519}{\line(-1,0){.0361271676}}
			\multiput(35.75,30.5)(-.0337338262,.0586876155){541}{\line(0,1){.0586876155}}
			\multiput(17.25,29)(.0336914063,.0356445313){512}{\line(0,1){.0356445313}}
			\multiput(34,62)(-.0337186898,-.0626204239){519}{\line(0,-1){.0626204239}}
			\put(14.25,25.5){\makebox(0,0)[cc]{$w_1$}}
			\put(35.75,26){\makebox(0,0)[cc]{$w_2$}}
			\put(14,62.25){\makebox(0,0)[cc]{$y_1$}}
			\put(35.75,64.75){\makebox(0,0)[cc]{$y_2$}}
			\put(18.25,76.5){\line(1,0){15.5}}
			\qbezier(53.75,77)(92.375,-8)(16.5,30)
		\end{picture}

		\end{center}
\vskip -2 cm
\caption{The graph $D^{odd}_{3, 9}$.}
\label{g39}
\end{figure}

We next show that every graph $D^{odd}_{k, n}$ is $D(k)$.
\vskip 5 pt

\begin{lem}\label{oddk}
For an odd integer $k \geq 3$ and a natural number $n \geq 4k - 3$, if $G$ is the graph $D^{odd}_{k, n}$, then $G$ is a $D(k)$ graph of order $n$.
\end{lem}
\proof For $1 \leq i \leq \frac{k-1}{2}$, let $B_i= P_{2i-1} \cup P_{2i}$. It can be observed that coloring each independent set $P_{1}, ..., P_{k}$ with a unique color will be a $\chi_{d}$-coloring of $G$. Hence, $\chi(G)=\chi_{d}(G) = k$.
\vskip 5 pt

\indent We will show that $\gamma(G)=k$. Clearly $\{x_1, x_2, \cdots, x_{k-1}, x_{k}\}$ is a dominating set of $G$. By the minimality of $\gamma(G)$, we have that $\gamma(G) \leq |\{x_1, x_2, \cdots, x_{k-1}, x_{k}\}|=k$. We next show that $\gamma(G) \geq k$. Let $D$ be a smallest dominating set of $G$. Because $D$ dominates $\{y_{2i},z_{2i},y_{2i-1},z_{2i-1}\}$, we have the following claim without proof.
\vskip 5 pt

\begin{cl}\label{claim1}
For $1 \leq i \leq \frac{k-1}{2}$, we have that $|D\cap B_{i}|\geq 2$.
\end{cl}
\vskip 5 pt

\indent We first assume that $x_{k} \in D$. Thus
\begin{align*}
\gamma(G) = |D|&=\displaystyle \sum_{i=1}^{\frac{k-1}{2}} |D\cap B_{i}|+|D\cap\{x_{k}\}|\geq 2\left(\frac{k-1}{2}\right)+1=k
\end{align*}
\noindent and this completes the proof.
\vskip 5 pt

\indent Hence, we may assume that $x_{k} \notin D$. We need the following claim.
\vskip 5 pt

\begin{cl}\label{claim2}
For $1 \leq i \leq \frac{k-1}{2}$, if $x_{2i-1}$ or $x_{2i}$ belongs to $D$, then $|D\cap B_{i}|\geq 3$.
\end{cl}
\textbf{Proof of Claim \ref{claim2}}. Assume $x_{2i-1} \in D$. Clearly to dominate $\{y_{2i-1},z_{2i-1},w_{2i-1}, w_{2i}\}$, $|D\cap (B_{i}\setminus \{x_{2i-1}\})| \geq 2$ implying that $|D\cap B_{i}|\geq 3$. The case when $x_{2i} \in D$ is proved similarly. This proves Claim \ref{claim2}.
\smallqed
\vskip 5 pt

\indent By Claim \ref{claim2}, if there exists $i_{0} \in \{1,2, \cdots, \frac{k-1}{2}\}$ such that $x_{2i_0-1} \in D$ or $x_{2i_0} \in D$, then $|D\cap B_{i_0}|\geq 3$. By Claim \ref{claim1}, we have
\begin{align*}
\gamma(G) = |D| =\displaystyle \sum_{i \neq i_0} |D\cap B_{i}|+|D\cap B_{i_0}| \geq 2\left(\frac{k-1}{2}-1\right)+3=k
\end{align*}	
\noindent and this completes the proof.
\vskip 5 pt

\indent Thus, we assume $x_{2i-1}, x_{2i} \notin D$ for all $1 \leq i \leq \frac{k-1}{2}$. The next claim is proved by similar arguments as Claim \ref{claim2} under the assumption $x_{k}, x_{2i-1}, x_{2i} \notin D$.
\vskip 5 pt

\begin{cl}\label{claim3}
For $1 \leq i \leq \frac{k-1}{2}$, if $w_{2i-1}$ or $w_{2i}$ belongs to $D$, then $|D\cap B_{i}|\geq 3$. Further, if $U \cap D \neq \emptyset$, then $|D\cap B_{1}|\geq 3$.
\end{cl}
\textbf{Proof of Claim \ref{claim3}}. For $1 \leq i \leq \frac{k-1}{2}$, we assume $w_{2i-1} \in D$. Because $x_{2i-1}, x_{2i} \notin D$, to dominate $\{y_{2i-1},z_{2i-1}, x_{2i-1}, x_{2i}\}$, it follows that $|D\cap (B_{i}\setminus \{w_{2i-1}\})| \geq 2$ Thus, $|D\cap B_{i}|\geq 3$. The case when $w_{2i} \in D$ is proved similarly.
\vskip 5 pt

It remains to consider the case when $D \cap U \neq \emptyset$. Let $u \in U \cap D$. To dominate $\{x_{1}, x_{2}, y_{1}, z_{1}, w_{1}, w_{2}\}$, we have that $|D\cap (B_{1}\setminus \{u\})| \geq 2$ which implies that $|D\cap B_{1}|\geq 3$. This proves Claim \ref{claim3}.
\smallqed
\vskip 5 pt

\indent Recall the assumption that $x_{k}, x_{2i-1}, x_{2i} \notin D$ for all $1 \leq i \leq \frac{k-1}{2}$. Since $D$ dominates $x_{k}$, $D \cap (W \cup U) \neq \emptyset$. Hence, $u \in D \cap U$ or there exists $i_{0} \in \{1,2, \cdots, \frac{k-1}{2}\}$ such that $D \cap \{w_{2i_0-1}, w_{2i_0}\} \neq \emptyset$. In both cases, we have $|D \cap B_{j_{0}}| \geq 3$ for some $j_{0} \in \{1, i_{0}\}$. Thus,
\begin{align*}
\gamma(G) = |D| =\displaystyle \sum_{i \neq j_0} |D\cap B_{i}|+|D\cap B_{j_0}| \geq 2\left(\frac{k-1}{2}-1\right)+3=k
\end{align*}	
and this proves Lemma \ref{oddk}.
\qed
\vskip 5 pt

\noindent \textbf{The graph} $D^{even}_{k, n}$
\vskip 5 pt

\noindent For an even integer $k \geq 4$, we let
\vskip 5 pt

\indent $X=\{x_1, x_2, \cdots, x_{k} \}$,
\vskip 5 pt

\indent $Y=\{y_1, y_2, \cdots, y_{k} \}$,
\vskip 5 pt

\indent $Z=\{z_1, z_2, \cdots, z_{k} \}$
\vskip 5 pt

\noindent and, for integer $t \geq 1$, we let
\vskip 5 pt

\begin{equation*}
U = \left\{
\begin{array}{rl}
	\emptyset & \text{if } n = 3k,\\\\
	\{u_{1}, ..., u_{t}\} & \text{if } n = 3k + t.
\end{array} \right.
\end{equation*}

\noindent Further, for $2 \leq i \leq k$ we let $P_i=\{x_i, y_i, z_i\}$  and $P_{1} = \{x_1, y_1, z_1\} \cup U$ be independent sets. The graph $D^{even}_{k, n}$ is obtained from $P_1, ..., P_{k}$ by adding edges as follows.
\begin{enumerate}
  \item For $1 \leq i \leq \frac{k}{2}$, join every vertex in $P_{2i}$ to every vertex in $P_{2i - 1}$.
  \item Join edges so that $X$ is $k$-clique.
\end{enumerate}

\noindent It can be observed by (a) that the subgraph induced by $P_{2i-1} \cup P_{2i}$ is a complete bipartite graph whose both partite sets have at least $3$ vertices for all $1 \leq i \leq \frac{k}{2}$. An example of the graph $D^{even}_{4, 12}$. is illustrated by Figure \ref{g412}.

\begin{figure}[H]

	\begin{center}
	\unitlength 1mm 
	\linethickness{0.4pt}
	\ifx\plotpoint\undefined\newsavebox{\plotpoint}\fi 
	\begin{picture}(99.75,91.25)(0,0)
		\put(21.5,52){\circle*{2.915}}
		\put(45.25,52.75){\circle*{2.915}}
		\put(69.25,52.75){\circle*{2.915}}
		\put(97,53){\circle*{2.915}}
		\put(22,35.5){\circle*{2.915}}
		\put(45.75,36.25){\circle*{2.915}}
		\put(69.75,36.25){\circle*{2.915}}
		\put(97.5,36.5){\circle*{2.915}}
		\put(21.75,17.75){\circle*{2.915}}
		\put(45.5,18.5){\circle*{2.915}}
		\put(69.5,18.5){\circle*{2.915}}
		\put(97.25,18.75){\circle*{2.915}}
		\qbezier(21.25,53)(55.25,91.25)(98.25,53.5)
		\qbezier(21,53.25)(49.375,73.625)(69.25,53.5)
		\qbezier(44.5,53)(66.875,77.875)(96.75,53.25)
		\multiput(96.75,53.25)(-1.8666667,-.0333333){15}{\line(-1,0){1.8666667}}
		\multiput(68.75,52.75)(-.03125,-.03125){8}{\line(0,-1){.03125}}
		\multiput(69,52.5)(-3,.03125){8}{\line(-1,0){3}}
		\multiput(45.25,53.25)(-1.076087,-.0326087){23}{\line(-1,0){1.076087}}
		\multiput(21.75,51.75)(.0525442478,-.0337389381){452}{\line(1,0){.0525442478}}
		\multiput(21.5,52.75)(.03373579545,-.04794034091){704}{\line(0,-1){.04794034091}}
		\multiput(21.75,36)(.0452793834,.0337186898){519}{\line(1,0){.0452793834}}
		\multiput(20.75,33.75)(.05,.0333333){15}{\line(1,0){.05}}
		\multiput(21.25,17.75)(1.0217391,.0326087){23}{\line(1,0){1.0217391}}
		\multiput(22,18.5)(.0440074906,.0337078652){534}{\line(1,0){.0440074906}}
		\multiput(45.25,52.75)(-.03372739917,-.04728789986){719}{\line(0,-1){.04728789986}}
		\multiput(22.25,35.5)(.61184211,.03289474){38}{\line(1,0){.61184211}}
		\multiput(21.5,35.5)(.0487551867,-.0337136929){482}{\line(1,0){.0487551867}}
		\multiput(69.5,53.5)(.0570539419,-.0337136929){482}{\line(1,0){.0570539419}}
		\multiput(69.25,53.5)(.03370098039,-.04197303922){816}{\line(0,-1){.04197303922}}
		\multiput(70,36.75)(.0540816327,.0336734694){490}{\line(1,0){.0540816327}}
		\put(69.75,36.5){\line(1,0){27.75}}
		\multiput(70,35.75)(.0540816327,-.0336734694){490}{\line(1,0){.0540816327}}
		\multiput(69,18.5)(1.2282609,.0326087){23}{\line(1,0){1.2282609}}
		\multiput(69.25,19)(.0529026217,.0337078652){534}{\line(1,0){.0529026217}}
		\multiput(97,52.5)(-.0337181045,-.04040097205){823}{\line(0,-1){.04040097205}}
		\put(16,54.5){\makebox(0,0)[cc]{$x_1$}}
		\put(43.75,57){\makebox(0,0)[cc]{$x_2$}}
		\put(71.75,56){\makebox(0,0)[cc]{$x_3$}}
		\put(98,57.5){\makebox(0,0)[cc]{$x_4$}}
		\put(17,36.25){\makebox(0,0)[cc]{$y_1$}}
		\put(46.5,39){\makebox(0,0)[cc]{$y_2$}}
		\put(68.25,39.25){\makebox(0,0)[cc]{$y_3$}}
		\put(98,39.25){\makebox(0,0)[cc]{$y_4$}}
		\put(17.75,17.75){\makebox(0,0)[cc]{$z_1$}}
		\put(48.75,20){\makebox(0,0)[cc]{$z_2$}}
		\put(69.75,15){\makebox(0,0)[cc]{$z_3$}}
		\put(99.75,16){\makebox(0,0)[cc]{$z_4$}}
	\end{picture}
	
	\end{center}	

\vskip -1 cm
\caption{The graph $D^{even}_{4, 12}$.}
\label{g412}
\end{figure}

We next show that every graph $D^{even}_{k, n}$ is $D(k)$.
\vskip 5 pt

\begin{lem}\label{evenk}
For an even integer $k \geq 4$ and a natural number $n \geq 3k$, if $G$ is the graph $D^{even}_{k, n}$, then $G$ is a $D(k)$ graph of order $n$.
\end{lem}
\proof Now for $1 \leq i \leq \frac{k}{2}$, let $B_i= P_{2i-1} \cup P_{2i}$. It can be checked that coloring each independent set with a unique color will be a $\chi_{d}$-coloring of $G$. Hence $\chi(G)=\chi_{d}(G)=k$. We will show that $\gamma(G)=k$. Clearly $\{x_1, x_2, \cdots, x_{k}\}$ is a dominating set of $G$. By the minimality of $\gamma(G)$, we have that $\gamma(G) \leq |\{x_1, x_2, \cdots , x_{k}\}|=k$. To show that $\gamma(G) \geq k$, we let $D$ be a smallest dominating set of $G$. Clearly, for $1 \leq i \leq \frac{k}{2}$,
\begin{align*}
|D\cap B_{i}| \geq 2.
\end{align*}
\noindent This implies that
\begin{align*}
\gamma(G) = |D| =\displaystyle \sum^{k/2}_{i = 1} |D\cap B_{i}| \geq 2\frac{k}{2} = k.
\end{align*}
Therefore $\gamma(G) \geq k$ and this proves Lemma \ref{evenk}.\\
\qed
\vskip 5 pt

\indent By Lemmas \ref{oddk} and \ref{evenk}, we establish realizability of a $D(k)$ graph of order $n \geq 4k - 3$ for all $k \geq 3$.
\vskip 5 pt

\noindent Recall that
\vskip 5 pt

\indent $\mathcal{D}(k, n):$ the class of $D(k)$ graphs of order $n$.
\vskip 5 pt

\noindent Now, we are ready to prove Theorem \ref{real}

\noindent \textbf{Proof of Theorem \ref{real}.} Assume that $n = 4k - 3 + t$ for some $t \geq 0$. When $k$ is odd, Lemma \ref{oddk} gives that $D^{odd}_{k, n} \in \mathcal{D}(k, n)$. We next consider the case when $k$ is even. Since $k \geq 3$, it follows that $4k - 3 \geq 3k$. Thus, Lemma \ref{evenk} gives that $D^{even}_{k, n} \in \mathcal{D}(k, n)$. This proves Theorem \ref{real}.

\subsection{Proof of Theorem \ref{thm d3}}

\begin{lem}\label{lemd3}
Let $G$ be a $D(3)$ graph with $\chi_{d}$-coloring $C=(V_1, V_2, V_3)$. Then $C$ has exactly one singleton color class.
\end{lem}
\proof Let $x_i \in V_i (1\leq i \leq 3)$. Suppose $|V_1|=|V_2|=1$, then $\{x_1,x_2\}$ is the dominating set of $G$, a contradiction. Therefore $C$ has at most one singleton color class. Now suppose $|V_i|\geq 2$, for all $1\leq i \leq 3$. Then the set $\{x_1, x_2, x_3\}$ is a minimum dominating set of $G$. Then each $x_i$ has at least one external private neighbor, let $y_i$ be an external private neigbor of $x_i$. If two vertices of $S=\{y_1,y_2,y_3\}$ belong to the same color class say $y_1, y_3\in V_2$ and $y_2\in V_1$, then $x_3$ has no dominating color class. Therefore each $y_i$ belongs to distinct color classes. Let $y_1 \in V_2 ; y_2\in V_3$ and $y_3 \in V_1$, then $\{y_1, y_2\}$ is the dominating set of $G$ since $y_1$ dominates $V_1$, $y_2$ dominates $V_2$ and any vertex of $V_3$ dominates $V_1$ or $V_2$ (not both) and hence adjacent to $y_1$ or $y_2$, a contradiction. Therefore $C$ has exactly one singleton color class.
\qed

\vskip 5 pt

\indent Now, we are ready to prove Theorem \ref{thm d3}.
\vskip 5 pt

\noindent \textbf{Proof of Theorem \ref{thm d3}.} Let $G$ be a $D(3)$ graph with $\chi_{d}$-coloring $C=(V_1, V_2, V_3)$. By the Corollary 2, let $\{x_i:  x_i\in V_i\}$ be the total dominating set of $G$ and let $y_i$ be an external private neighbor of $x_i$. By Lemma \ref{lemd3}, $C$ has exactly one singleton color class, let $V_{3}=\{x_3\}$. Clearly $y_i$ dominates $V_i$. Let $y_1, y_3 \in V_2$ and $y_2 \in V_1$. Thus, $|V_2|\geq 3$. Since $x_1$ cannot dominate $V_2$, it follws that $x_1$ dominates $V_3$ and it may be adjacent to some (or all) the vertices of $V_{2}-\{y_3\}$. Similarly $y_3$ is adjacent to $x_3$ and some (or all) the vertices of $V_{1}-\{x_1\}$. Now suppose $|V_1|=2$, then $\{x_1, y_2\}$ is the dominating set of $G$, a contradiction. Therefore $|V_1|\geq 3$. Now by the Theorem 1, any vertex of $V_{1}-\{x_1,y_2\}$ dominates $V_2$ or $V_{3}$ (not both) and any vertex of $V_{2}-\{y_3, y_1\}$ dominates $V_1$ or $V_{3}$(not both). If $x_3$ has exactly one non-neighbor in $V_2$ (i.e $y_{1}$). Then $\{x_3, y_1\}$ is dominating set of $G$, a contradiction. Similarly $x_3$ cannot have exactly one non-neighbor in $V_1$. Now if there exists $x \in V_1$ and $y \in V_2$ where $xy \notin E(G)$ such that $x$ has exactly one non-neighbor in $V_2$ (i.e $y$) and $y$ has exactly one non-neighbor in $V_1$ (i.e $x$), then $\{x,y\}$ is a dominating set of $G$ contradicting the minimality of $\gamma(G)$.
\vskip 5 pt

\indent Conversely let $G$ be constructed from three independent sets $V_1, V_2, V_3$ with $|V_1|\geq 3, |V_2|\geq 3, |V_3|=1$ satisfying the construction of the class $\mathcal{D}(3)$. Because $G$ contains a 5-cycle $y_2y_1x_1x_3y_3y_2$, $G$ is not bipartite. This implies that $\chi(G)=3$. Next, coloring the vertices of $V_1,V_2,V_3$ with different colors gives a $\chi_{d}$-coloring of $G$. Therefore $\chi_{d}(G)=3$.
\vskip 5 pt

\indent We will show that $\gamma(G)=3$. Suppose to the contrary that there exists a set $\{x,y\}$ which is a dominating set of $G$. Clearly, $x,y$ cannot belong the same color class $V_1$ or $V_2$ since $|V_1|, |V_{2}| \geq 3$. Assume that $x \in V_i$ for some $i \in \{1, 2\}$ and $y \in V_3$. By the construction, $y = x_{3}$. Since $\{x, y\}$ is a dominating set of $G$, it follows that $x_{3}$ dominates $V_{i} \setminus \{x\}$. This violates Property (d) of the construction that $x_{3}$ has at least two non-neighbors in each of $V_{1}$ and $V_{2}$.
\vskip 5 pt

\indent Hence, we assume that $x \in V_1$ and $y \in V_2$. If $x$ dominates $V_2$ and $y$ dominates $V_1$, then Property (d) of the construction implies that both $x$ and $y$ are not adjacent to $x_{3}$. Thus, $x_3$ is not dominated by $\{x, y\}$, a contradiction. Thus, renaming vertices if necessary, we assume that $y$ does not dominate $V_1$. If $xy \in E(G)$, then there is a vertex $v \in V_{1} \setminus \{x\}$ such that both $x$ and $y$ are not adjacent to $v$. This contradicts $\{x, y\}$ is a dominating set. So, we assume that $xy \notin E(G)$. Thus, Property (e) of the construction yields that there are at least two non-neighbors of $y$ in $V_{1}$ or there are at least non-neighbors of $x$ in $V_{2}$. This contradicts $\{x, y\}$ is a dominating set. Therefore $\gamma(G) = 3$. This establishes the characterization of $D(3)$ graphs.
\vskip 5 pt

\indent Finally, to show that there is no planar $D(3)$ graphs, we let $G$ be a planar graph. Suppose to the contrary that $G$ is $D(3)$. Thus $G \in \mathcal{D}(3)$. Let $C=\{V_1, V_2, V_3\}$ ($|V_1|\geq 3, |V_2|\geq 3, |V_3|=1$) be a $\chi_{d}$-coloring of $G$ such that$\{x_i:  x_i\in V_i\}$ be the total dominating set of $G$ (by the Corollary 2) and let $y_i$ be the external private neighbor of $x_i$. Let $z_1$ and $z_2$ be adjacent to $x_3$. Let $y_1, y_3 \in V_2$ and $y_2 \in V_1$. Since $G \in \mathcal{D}(3)$, $x_3$ has at least two non- neighbors in $V_1$ (say $u_1, v_1$) and $x_3$ has at least two non- neighbors in $V_2$ (say $u_2, v_2$). Then $G[\{u_1, v_1, z_1\} \cup \{u_2, v_2, z_2\}]$ is a minor of $K_{3,3}$ since there is path $z_1  x_3  z_2$ which contradicts the fact that $G$ is planar. This proves Theorem \ref{thm d3}.

\section{Proof of Theorem \ref{thm 2}}

Let $G$ be a graph with $\chi_{d}$-coloring $C=(V_1, V_2,\cdots, V_k)$ with $k \geq 3$. For $1 \leq i, j ,l \leq k$ such that $i\neq j \neq l$, a \emph{chain} $V_i \rightarrow V_j \rightarrow V_l$ means there exist a vertex $x_i \in V_i$ such that $x_i$ dominates $V_j$, a vertex $x_j \in V_j$ such that $x_j$ dominates $V_l$ and a vertex $x_l \in V_l$ such that $x_l$ dominates $V_i$.
\vskip 5 pt

\begin{lem}\label{lem5}
Let $G$ be a graph with $\chi_{d}$-coloring $C=(V_1, V_2, V_3, V_4)$ such that $|V_i| \geq 2$ for all $1 \leq i \leq 4$. If there exists a chain $V_i \rightarrow V_j \rightarrow V_k$, then $G$ is not a $D(4)$ graph.
\end{lem}
\proof If there exists a chain $V_i \rightarrow V_j \rightarrow V_k$, then the set $D=\{x_i, x_j, x_k : x_i \in V_j, x_j \in V_j, x_k \in V_k\}$ dominates the color classes $V_i, V_j, V_k$. By Theorem \ref{thm 1} any vertex in $V_l$ dominates exactly one color class $V_i, V_j$ or $V_k$ where $l \neq i, j, k$. Thus, every vertex of $V_l$ has a neighbor in $D$. Therefore $\gamma(G) \leq 3$ and hence $G$ is not a $D(4)$ graph.
\qed
\vskip 5 pt

\noindent \textbf{Proof of Theorem \ref{thm 2}.} Suppose to the contrary that $G$ is $D(4)$. Thus, $G$ has a $\chi_{d}$-coloring $C=(V_1, V_2, V_3, V_4)$ of $G$. Renaming the color classes if necessary, we assume that $|V_{1}| \geq |V_{2}| \geq |V_{3}| \geq |V_{4}|$. Because $\gamma(G) = 4$, at most two color classes are singleton. Thus, $|V_{1}|, |V_{2}| \geq 2$. We distinguish $7$ cases due to the number of vertices in each color class.
\vskip 5 pt

\noindent \textbf{Case 1:} $|V_3|=1$ and $|V_4|=1$.\\
\indent Let $V_3 = \{x_3\}, V_4 = \{x_4\}$ and let $U_i=V_{i}-(N(x_3) \cup N(x_4) )$ for $i= 1, 2$. Since $\chi(G)=4$, $x_3 x_4 \in E(G)$ and both $x_3, x_4$ have at least one neighbor in $V_1$ and $V_2$. Now if $|U_1| \geq 1$ and $|U_2| \geq 1$ (say $x_1 \in U_1$ and $x_2 \in U_2$), then $x_1$ dominates $V_2$ and $x_2$ dominates $V_1$ implying that $\{x_1, x_2, x_3\}$ is a dominating set of $G$ of cardinality 3, a contradiction. So let us assume that $|U_1| \geq 1$ and $|U_2| = 0$. Then every vertex of $U_1$ dominates $V_2$. Let $y_2$ be the neighbor of $x_3$ in $V_2$. Then the set $\{y_2, x_3, x_4\}$ is a dominating set of $G$. If $|U_1|= |U_2| = 0$, then $\{x_3, x_4\}$ is a dominating set of $G$. Both contradicts $\gamma(G) = 4$. Therefore, Case 1 cannot occur.
\vskip 5 pt

\noindent \textbf{Case 2:} $|V_4|=1$ and $|V_i|\geq 2$ for $i=1,2,3$.\\
\indent Let $V_4 = \{x_4\}$. Let $U_i=V_{i}- N(x_4)$ for $i= 1, 2, 3$. If $|U_i| \geq 3$ for all $1 \leq i \leq 3$, then by pigeonhole principle at least two vertices in $U_i$ dominate $U_j$. Let two vertices of $U_1$ (say $u_1, u_2$) dominate $U_2$. If every vertex in $U_{3}$ dominates $U_{1}$, then $G[U_{1} \cup U_{3}]$ contains $K_{3, 3}$ which is a contradiction. Thus, at least one vertex of $U_3$ (say $u_3$) dominate $U_2$. Then $G[U_{2} \cup \{u_1,u_2,u_3\}]$ is $K_{3,3}$ contradicting that $G$ is planar.
\vskip 5 pt

\indent Let $x_i$ be the neighbor of $x_4$ in $V_i$, for $1 \leq i \leq 3$. Thus, $\{x_i: 1\leq i \leq 4\}$ is a total dominating set of $G$ implying that every $x_i$ has an external private neighbor $y_i$ with respect to $\{x_{1}, x_{2}, x_{3}, x_{4}\}$. So, $y_i$ dominates $V_i$.
\vskip 5 pt

\indent We first assume that $|U_1| \leq 2$ and $|U_i| \geq 3$ for $i=2,3$. We consider the case when $|U_1|=1$. Let $U_1 = \{u_1\}$. If $u_1$ dominates $V_2$, then $\{u_1, y_3, x_4\}$ is a dominating set of $G$. If $u_1$ dominates $V_3$, then $\{u_1, y_2, x_4\}$ is a dominating set of $G$. Both contradicts $\gamma(G) = 4$. We consider the case when $|U_1|=2$. Let $U_1 = \{u_1, v_1\}$. If $u_1$ dominates $V_2$ and $v_1$ dominates $V_3$, then $\{u_1, v_1, x_4\}$ is a dominating set of $G$, a contradiction. Therefore both $u_1, v_1$ dominate the same color class. Assume that each of $u_1, v_1$ dominates $V_2$. If every vertex of $U_3$ dominate $V_1$, then $\{u_1, v_1, x_4\}$ is a dominating set of $G$, a contradiction. So there exists a vertex say $u_3 \in U_3$ such that $u_3$ dominates $V_2$. then $G[U_{2} \cup \{u_1,v_1,u_3\}]$ is $K_{3,3}$ contradicting $G$ is planar.
\vskip 5 pt

\indent Now let us assume that $|U_1| \leq 2, |U_2|\leq 2$. We prove that this case does not occur whether $|U_{3}| \geq 3$ or $|U_{3}| \leq 2$. If $|U_1|=|U_2|=1$. Let $u_1 \in U_1$ and $u_2 \in U_2$. If $u_1 u_2 \in E(G)$, then $\{u_1, y_3, x_4\}$ is a dominating set of $G$. If $u_1 u_2 \notin E(G)$, then both $u_1, u_2$ dominate $V_3$ implying that $\{u_1, u_2, x_4\}$ is a dominating set of $G$. Both contradict $\gamma(G) = 4$. So let us assume that $|U_1|=1$ and $|U_2|=2$. If $u_1$ dominates $V_2$, then $\{u_1, y_3, x_4\}$ is a dominating set of $G$. If $u_1$ dominates $V_3$, then $\{u_1, y_2, x_4\}$ is a dominating set of $G$. Both contradict $\gamma(G) = 4$. We prove that $|U_1|=2$ and $|U_2|=1$ cannot occur by similar arguments. Hence, we assume that $|U_{1}| = |U_{2}| = 2$. Thus, $|V_{2}| \geq 3$ as there is $x_{2} \in V_{2} \setminus U_{2}$. Let $U_{1} = \{u_{1}, v_{1}\}$. We consider when $u_{1}$ and $v_{1}$ dominate different color classes. Without loss of generality, we assume that $u_{1}$ dominates $V_{2}$ and $v_{1}$ dominates $V_{3}$. So, $\{x_{4}, u_{1}, v_{1}\}$ is a dominating set of $G$ contradicting $\gamma(G) = 4$. Thus, each of $u_{1}, v_{1}$ dominates $V_{j}$ for some $j \in \{2, 3\}$. It can be proved that the case when $|U_{3}| \leq 1$ cannot occur by similar arguments to the case when $|U_1|=1$ and $|U_2|=2$. Thus, we can assume that $|U_{3}| \geq 2$ which yields that $|V_{3}| \geq 3$. Since $G$ is $K_{3, 3}$-free, every vertex $u$ in $U_{5 - j}$ dominates $V_{1}$ implying that $\{x_{4}, u_{1}, v_{1}\}$ is a dominating set of $G$, a contradiction. Therefore, Case 2 cannot occur.
\vskip 5 pt

\noindent \textbf{Case 3:} $|V_i|=2$ for all $1 \leq i \leq 4$.\\
\indent Let $x_i, y_i \in V_i$. Since $\chi(G)=4$, there exists a vertex in each color that has at least one neighbor in every other color class. Let $x_1$ be adjacent to $x_2, x_3, x_4$. Assume that each of $x_{1}, y_{1}$ dominates different color classes. Renaming vertices if necessary, we let $x_1$ dominates $V_2$ and $y_1$ dominates $V_3$. So, $\{x_1, y_1, y_4\}$ is a dominating set of $G$ contradicting the minimality of $\gamma(G)$. Therefore, each of $x_1, y_1$ dominates the same color class, $V_2$ say. If each of $x_3, y_3$ dominates $V_i$ for some $i \in \{1, 2, 4\}$, then $\{x_1, y_2, y_4\}$ is a dominating set of $G$ contradicting the minimality of $\gamma(G)$. So let us assume that $x_3$ dominates $V_j$ and $y_3$ dominates $V_{j'}$ where $j, j' \in \{1, 2, 4\}$and $j \neq j'$. In all the cases, it can be checked that $\{x_1, y_2, y_4\}$ is a dominating set of $G$ contradicting the minimality of $\gamma(G)$. Therefore, Case 3 cannot occur.
\vskip 5 pt

\noindent \textbf{Case 4:} $|V_1| \geq 3$ and  $|V_i|=2$ for $2 \leq i \leq 4$.\\
\indent Let $x_i, y_i \in V_i$ for $1 \leq i \leq 4$. Let $x_1$ be adjacent to $x_2, x_3, x_4$. By Theorem \ref{thm 1}, the class $V_{1}$ is dominated by a vertex. Renaming vertices if necessary, we let $x_{2}$ or $y_{2}$ dominate $V_{1}$.
\vskip 5 pt

\indent We first assume that both of $x_2, y_2$ dominate $V_1$. Thus, no vertex of $V_3$ or $V_4$ dominate $V_1$ otherwise $G$ would contain a $K_{3,3}$. We may assume that $x_{3}$ dominates $V_{j}$ and $y_{3}$ dominates $V_{j'}$ where $j, j' \in \{2, 4\}$. We have that $\{x_1, x_2, y_4\}$ is a dominating set of $G$ contradicting the minimality of $\gamma(G)$.
\vskip 5 pt

\indent We next assume that $x_2$ dominates $V_1$ but $y_2$ dominates a different color class, $V_3$ say. Since $\{x_1,x_2,x_3,x_4\}$ is a total dominating set of $G$, let $z_i$ be an external private neighbor of $x_i$, for $1 \leq i \leq 4$. Then $\{x_2,y_2,z_4\}$ is a dominating set of $G$ contradicting the minimality of $\gamma(G)$. Case 4 cannot occur.
\vskip 5 pt

\noindent \textbf{Case 5:} Let $|V_1| \geq 3, |V_2| \geq 3$ and  $|V_i|=2$ for $i=3,4$.\\
\indent  Let $x_i, y_i \in V_i$ for $1 \leq i \leq 4$. Suppose that each of $x_1, y_1$ dominates $V_2$. Then no vertex of $V_3$ or $V_4$ dominate $V_2$, otherwise $G$ would contain a $K_{3,3}$. Similarly, every vertex of $V_1 - \{x_1, y_1\}$ dominate $V_3$ or $V_4$. If $y_3$ dominates $V_4$, then $\{x_1, x_3, y_4\}$ is a dominating set of $G$. If $y_3$ dominates $V_1$, then $\{x_1, y_3, y_4\}$ is a dominating set of $G$. Hence, we can assume that there is only vertex $x_1$ in $V_{1}$ that dominates $V_2$. So, every vertex of $V_1 - \{x_1\}$ dominates $V_3$ or $V_4$. Because $G$ does not contain $K_{3, 3}$, it follows that, for each $i \in \{3, 4\}$, there exits $w_{i} \in V_{i}$ which does not dominate $V_{2}$. Thus, $w_{i}$ dominates $V_{7 - i}$ or $V_{1}$. Let $\{u_{i}\} = V_{i} - \{w_{i}\}$. Clearly, $\{x_{1}, u_{3}, u_{4}\}$ is a dominating set of $G$, a contradiction.
\vskip 5 pt

Now let us assume that no vertex of $V_1$ dominates $V_2$. Then every vertex of $V_1$ dominates of $V_3$ or $V_4$ and every vertex of $V_2$ dominates $V_1$ or $V_3$ or $V_4$. Then $\{x_1, y_3, y_4\}$ is a dominating set of $G$. Therefore, Case 5 cannot occur.
\vskip 5 pt

\noindent \textbf{Case 6:} Let $|V_1| \geq 3$, $|V_2| \geq 3$, $|V_3|\geq 3$ and  $|V_4|=2$ .\\
\indent  $x_i, y_i \in V_i$ for $1 \leq i \leq 4$. Let $x_1$ be adjacent to $x_2, x_3, x_4$. Suppose each of $x_4, y_4$ dominates the same color class, $V_3$ say. Then no vertex of $V_1$ and $V_2$ dominate $V_3$, otherwise there would be a $K_{3,3}$. Let $S_1$ be the set of vertices in $V_1$ that dominate $V_2$ and $S_2$ be the set of vertices in $V_2$ that dominate $V_1$. It is clear that $1 \leq |S_1|\leq 2$ and $1 \leq |S_2|\leq 2$, otherwise there would be a $K_{3,3}$. every vertex of $V_1 - S_1$ and $V_2 - S_2$ dominates $V_4$. If $|S_1|=1$, then let $s_1 \in S_1$ and $s_2 \in S_2$. Thus, $\{s_1, x_3, y_4\}$ is a dominating set of $G$. Similarly if $|S_2|=1$, we arrive at the contradiction. Let us assume that $|S_1|=|S_2|=2$. Let $S_i = \{s_i, r_i\}$, for $i=1,2$. Let $v_i \in V_i -S_i$. Now the graph $G[\{s_1, r_1, v_1\} \cup \{s_2, r_2, v_2\}]$ is a minor of $K_{3,3}$ since there is path $v_1 x_4 x_3 y_4 v_2$. Now if $x_4$ dominates $V_3$ and $y_{4}$ dominates $V_2$, then Theorem \ref{thm 1} yields that there exists a vertex $z$ that dominates $V_{1}$. Thus, $\{x_4, z, y_4\}$ is a dominating set of $G$ contradicting the minimality of $\gamma(G)$. Therefore, Case 6 cannot occur.
\vskip 5 pt

\noindent \textbf{Case 7:} Let $|V_i| \geq 3$ for all $1 \leq i \leq 4$.\\
\indent By Theorem \ref{thm 1} every vertex dominates exactly one color class. We shall prove the following claim.
\vskip 5 pt

\begin{cl}\label{claim4}
For each $1 \leq i \leq 4$, every vertex in $V_i$ dominates the same color class.
\end{cl}
\noindent \textbf{Proof of Claim \ref{claim4}}. Suppose to the contrary that there exist two vertices $x, y \in V_2$ such that $x$ dominates $V_1$, $y$ dominates $V_3$. Thus, $y$ is not adjacent to some vertex in $V_{1}$, $a_{1}$ say. By Theorem \ref{thm 1}, $a_{1}$ dominates $V_3$ or $V_4$. We distinguish 2 cases.
\vskip 5 pt

\noindent \textbf{Subcase (i):} \emph{$a_1$ dominates $V_3$.}\\
\indent If there exists a vertex in $V_3$ that dominates $V_2$, then $G$ has a chain $V_3\rightarrow V_2\rightarrow V_1$ contradicting Lemma \ref{lem5}. Hence, every vertex in $V_3$ does not dominate $V_2$. Thus, by Theorem \ref{thm 1}, every vertex in $V_{3}$ dominates either $V_{1}$ or $V_{4}$.
\vskip 5 pt

\indent So, we suppose that there is a vertex of $V_3$ that dominates $V_1$. If there exists a vertex in $V_1$ that dominates $V_2$, then $G$ has a chain $V_1\rightarrow V_2\rightarrow V_3$, contradicting Lemma \ref{lem5}. So every vertex in $V_1$ does not dominate $V_2$. Since $G$ is planar, the subgraph induced on $V_1$ and $V_3$ is not complete bipartite. Thus, there is a vertex $b \in V_{3}$ which is not adjacent to a vertex in $V_{1}$. By Theorem \ref{thm 1}, $b$ dominates $V_4$. If there exists a vertex in $V_4$ that dominates $V_2$, then there exists a chain $V_4\rightarrow V_2\rightarrow V_3$. If there exists a vertex in $V_4$ that dominates $V_1$, then there exists a chain $V_4\rightarrow V_1\rightarrow V_3$. Both contradict Lemma \ref{lem5}. Thus, by Theorem \ref{thm 1}, every vertex in $V_{4}$ is adjacent to every vertex in $V_{3}$ yielding that $G[V_{3} \cup V_{4}]$ has $K_{3, 3}$ as a subgraph. This contradicts $G$ is a planar. Subcase (i) cannot occur.\\

\noindent \textbf{Subcase (ii):} \emph{$a_1$ does not dominates $V_3$.}\\
\indent Thus, by Theorem \ref{thm 1}$, a_1$ dominates $V_{4}$. Note that no vertex of $V_1$ dominates $V_3$ (and vice-verse) and no vertex of $V_2$ dominates $V_4$ (and vice-versa). Otherwise it would reduce to Subcase (i) or $G$ would have a chain. So every vertex in $V_{1}$ dominates either $V_{2}$ or $V_{4}$. Since $|V_{2}| \geq 3$, by the pigeonhole principle, there are at least two vertices in $V_1$, say $u, v$ such that both $u, v$ either dominate $V_2$ or dominate $V_4$.
\vskip 5 pt

\indent We first consider the case when $u, v$ dominate $V_2$. By Theorem \ref{thm 1}, every vertex in $V_{3}$ dominates either $V_{2}$ or $V_{4}$. If every vertex in $V_{3}$ dominates $V_{4}$, then $G[V_{3} \cup V_{4}]$ contains $K_{3, 3}$ contradicting $G$ is planar. Thus, there is at least one vertex in $V_3$, say $w_1$ that dominates $V_2$. Thus, $G[\{u, v, w_1\} \cup V_2]$ is $K_{3,3}$ contradicting $G$ is planar.
\vskip 5 pt

\indent Hence, we consider the case when $u, v$ dominate $V_4$. Similarly, every vertex in $V_{3}$ dominates either $V_{2}$ or $V_{4}$. By planarity of $G$, there is at least one vertex in $V_3$, say $w_2$ that dominates $V_4$. Thus, $G[\{u, v, w_2\} \cup V_4]$ contains $K_{3,3}$ contradicting $G$ is planar. Subcase (ii) cannot occur and this proves Claim \ref{claim4}.
\vskip 5 pt

\indent We suppose by Claim \ref{claim4} that every vertex in $V_{i}$ dominates $V_{j}$. Hence, $G[V_{i} \cup V_{j}]$ is $K_{3, 3}$ contradicting $G$ is planar. Case 7 cannot occur. Therefore $G$ is not $D(4)$ and this proves Theorem \ref{thm 2}.

\section{Open problems}
1. For all $k\geq 3$, we have constructed $D^{odd}_{k, n}$ for all $n \geq 4k-3$ and $D^{even}_{k, n}$ for all $n \geq 3k$. We believe that these are smallest possible orders of $D(k)$ graphs in each parity. We conjecture that:
\vskip 5 pt

\begin{conj}
For all $k\geq 3$, there exists no $D^{odd}_{k, n}$ for all $n < 4k-3$ and $D^{even}_{k, n}$ for all $n < 3k$.
\end{conj}
\vskip 5pt

\indent By Theorem \ref{thm 5}, it would be interesting to characterize all $D(k)$ non-planar graphs. When $k = 2$, the $D(2)$ graphs which are non-planar are $K_{p, q}$ for $p \geq q \geq 3$. When $k = 3$, by Theorem \ref{thm d3}, a non-planar graph $G$ is $D(3)$ if and only if $G \in \mathcal{D}(3)$. We post an open problem that:

\begin{problem}
Characterize all $D(k)$ non-planar graphs for all $k \geq 4$.
\end{problem}

\section*{Acknowledgements}
The first author is thankful to UGC, New Delhi, for UGC-JRF, under which this work has been done.
The first and third authors are thankful to University Grants Commission (UGC), India for financial support under the grant UGC-SAP-DRS-II, NO.F.510/12/DRS-II/2018(SAP-I) dated: 9$^{th}$ April 2018.

\end{document}